\documentclass[a4paper,fleqn]{cas_sc}

\usepackage[numbers]{natbib}
\newtheorem{definition}{Definition}%
\newtheorem{proposition}{Proposition}
\newtheorem{corollary}{Corollary}
\newtheorem{theorem}{Theorem}
\newtheorem{lemma}{Lemma}
\newdefinition{rmk}{Remark}
\newproof{pf}{Proof}
\newproof{pot}{Proof of Theorem \ref{thm}}
\usepackage{graphicx}%
\usepackage{multirow}%
\usepackage{amsmath,amssymb,amsfonts}%
\usepackage{algorithm}
\usepackage{algpseudocode}
\usepackage{tikz}  
\begin{document}
\let\WriteBookmarks\relax
\def\floatpagepagefraction{1}
\def\textpagefraction{.001}

\shorttitle{}    

\shortauthors{}  

\title [mode = title]{Scalable inference of large-scale random kronecker graphs via tensor decomposition and Einstein summation}  




\author[1]{Sanaa Khobizy}


\fnmark[b]



\affiliation[1]{
    organization={Laboratory LAMAI, University Cadi Ayyad},
    city={Marrakech}, 
    country={Morocco}
}

\affiliation[2]{
    organization={Laboratory L.M.P.A, University Littoral, Côte d’Opale}, 
    city={Calais}, 
    country={France}
}




\ead{khobizy.sanaa23@gmail.com}





\begin{abstract}
In this paper, we extend the analysis of random Kronecker graphs to multi-dimensional networks represented as tensors, enabling a more detailed and nuanced understanding of complex network structures. We decompose the adjacency tensor of such networks into two components: a low-rank signal tensor that captures the essential network structure and a zero-mean noise tensor that accounts for random variations. Building on recent advancements in tensor decomposition and random tensor theory, we introduce a generalized denoise-and-solve framework that leverages the Einstein summation convention for efficient tensor operations. This approach significantly reduces computational complexity while demonstrating strong performance in network inference tasks, providing a scalable and efficient solution for analyzing large-scale, multi-dimensional networks.
\end{abstract}




\maketitle

\section{Introduction}\label{sec1}
In today's increasingly interconnected world, the amount of data generated by large-scale systems is growing rapidly\cite{myers2014information,newman2002random,pavlopoulos2011using,koutrouli2020guide}. These systems include social networks (e.g., Facebook, Twitter) \cite{barabasi1999emergence}, transportation networks (e.g., traffic management systems, ride-sharing services) \cite{newman2003structure}, and biological networks (e.g., protein-protein interaction networks, neural networks) \cite{jure2014snap}. The relationships between the components of these systems can be complex, and their understanding requires sophisticated analytical methods. One of the most effective ways to represent and analyze these systems is through graph models, which capture the relationships between entities as edges between nodes. Graphs are ideal for capturing the structure of such systems because they can represent entities as nodes (vertices) and their interactions or relationships as edges (links). Probabilistic graph models are especially useful for large-scale networks because they simplify the analysis by reducing the complexity of the system to only a few key parameters. These models rely on probabilistic principles to determine the existence of edges between nodes, making them capable of handling the randomness and uncertainty that are common in real-world networks. One of the foundational random graph models is the Erdős–Rényi graph, where each edge between two nodes is assigned a fixed probability, and the edges are formed randomly. This model, while foundational, does not capture many of the complexities seen in real-world networks, such as the presence of community structures (groups of nodes that are more densely connected to each other than to nodes outside the community) or the small-world phenomenon (a network where most nodes are not directly connected but can be reached through a small number of hops). To address these limitations, more sophisticated models like the stochastic block model and the Watts–Strogatz model have been developed \cite{karrer2011stochastic,watts1998collective,abbe2018community}. The stochastic block model is effective for capturing community structure in networks by grouping nodes into blocks (or communities), where nodes within the same block have a higher probability of being connected. The Watts–Strogatz model, on the other hand, focuses on small-world properties, where most nodes are not directly connected but can reach each other through a small number of intermediate nodes, creating a high level of clustering \cite{serrano2003topology}.

A more recent development in the analysis of complex networks is the generalized random Kronecker graph model. The strength of this model lies in its ability to generate large and complex graphs using relatively simple and efficient mathematical operations. The model relies on the Kronecker product, a mathematical operation that combines two smaller matrices (or tensors) to form a larger one \cite{leskovec2010kronecker,leskovec2008statistical,liao2023analysis}. By repeatedly applying this operation, the Kronecker graph model can generate large probability tensors that represent the structure of large networks. The primary advantage of the Kronecker graph model is that it can generate large, synthetic graphs that closely resemble real-world networks, making it highly valuable for simulation and testing. Researchers can use these synthetic graphs to evaluate algorithms, test hypotheses, or perform simulations before applying them to actual data. Moreover, the Kronecker model can be used for graph classification, where the model's parameters (the Kronecker initiator tensor) are estimated based on real network data. This allows for the identification and classification of network types, making the model useful in applications like social network analysis, transportation networks, and biological network modeling. While the generalized Kronecker graph model has proven effective in modeling real-world networks, there are several challenges associated with it. One of the primary issues is computational efficiency. The estimation of graph parameters (the values in the Kronecker initiator tensor) for large networks can be computationally expensive \cite{errica2019fair,todorovic2012human,fan2022big}. As the size of the network grows, the process of generating and analyzing these large graphs becomes increasingly resource-intensive, which can limit the model's scalability. Additionally, the Kronecker graph model tends to produce isomorphic graphs (graphs that are structurally identical), which can complicate the process of matching corresponding nodes across different graphs. This can make the graph inference process more difficult and computationally expensive. For example, when performing tasks like link prediction or graph classification, finding the correct correspondence between nodes in different graphs becomes more complex, requiring advanced techniques and computational resources. As networks grow larger and more complex, traditional graph models like the Erdős–Rényi and even Kronecker graphs must adapt to handle new challenges such as dynamic networks, where the structure changes over time, and multilayer networks, where nodes have multiple types of relationships with other nodes (e.g., Facebook users interacting through both "friend" and "group membership" edges). These new dimensions further complicate the graph analysis process, requiring enhanced models and computational techniques to effectively handle and analyze these evolving systems.
\section{Backgrounds and notations}\label{sec2}
In this section, we offer a brief summary of fundamental concepts and introduce important terms that will be used throughout the rest of the paper. A tensor is a multidimensional data structure, with its order determined by the number of indices, commonly known as modes or dimensions. Tensors generalize the notions of scalars, vectors, and matrices to higher-dimensional spaces. In particular, a scalar is a tensor of order zero, a vector is a tensor of order one, and a matrix is a tensor of order two. In this paper, lowercase letters (e.g. $a$) denote vectors, uppercase letters (e.g., $A$) represent matrices, and calligraphic letters (e.g., $\mathcal{A}$) are used for higher-order tensors. We start by reviewing the definition of the n-mode product between a tensor and a matrix; further information is available in \cite{kolda2009tensor,de2008decompositions}.
\begin{definition}
            The $n$-mode product of the tensor $\mathcal{A} = \left[a_{i_1 i_2 \ldots i_n}\right] \in \mathbb{R}^{I_1 \times \ldots \times I_N}$ and the matrix $U = \left[u_{j i_n}\right] \in \mathbb{R}^{J \times I_n}$, denoted as $\mathcal{A} \times_n U$, results in a tensor of order $I_1 \times \ldots \times I_{n-1} \times J \times I_{n+1} \times \ldots \times I_N$. The entries of this tensor are given by:

$$
\left(\mathcal{A} \times_n U\right)_{i_1 i_2 \ldots i_{n-1} j i_{n+1} \ldots i_N} = \sum_{i_n=1}^{I_n} a_{i_1 i_2 \ldots i_N} u_{j i_n}.
$$
   			\end{definition}
            \begin{definition}
	The n-mode vector product of the $ N $ order tensor  $ \mathcal{A}$ with a vector $ v \in \mathbb{R}^{I_n} $ is an $ N-1 $ order tensor of size $ I_1\times I_2\times\cdots\times I_{n-1}\times
	I_{n+1} \times \cdots\times I_N $ denoted by $ \mathcal{X }\bar{\times}_n v $ and given by:
	$$
	( \mathcal{A}\bar{\times}_n v)_{i_1 \ldots i_{n-1} i_{n+1} \ldots i_N }= \sum_{ i_n=1}^{I_n} \mathcal{A}(i_1,i_2, \ldots ,i_N) v_{i_n}.
	$$
\end{definition}
\begin{proposition}\cite{lee2014fundamental}
Let $\mathcal{A}\in \mathbb{R}^{R_1\times R_2\times...\times R_N }$ be an N-th order tensor. Then
\begin{enumerate}
    \item $\mathcal{A}\times_mB\times_nC=\mathcal{A}\times_nC\times_mB.$
    \item $\mathcal{A}\times_nB\times_nC=\mathcal{A}\times_nCB.$
    \item Moreover, If $B$ has full column rank, then
    $$ \mathcal{C}=\mathcal{A}\times_n B\Rightarrow\mathcal{A}=\mathcal{C}\times_nB^{+},$$
    where $B^{+}$ is the Moore-Penrose pseudoinverse of $B$.
    
    \item If $B\in\mathbb{R}^{R\times I_n}$, then
    $$ \mathcal{C}=\mathcal{A}\times_n B\Leftrightarrow C_{(n)}=BA_{(n)}.$$
\end{enumerate}
\end{proposition}
      \noindent
  We now revisit the definition and key properties of the tensor Einstein product, which serves as a generalization of matrix multiplication to higher-dimensional tensors. The Einstein product provides a framework for performing multiplication operations between tensors, extending the familiar concept of matrix multiplication to handle multidimensional arrays. This operation is crucial in various fields, including machine learning, signal processing, and tensor decomposition, as it allows for the efficient manipulation and analysis of complex data structures. For a comprehensive understanding of the tensor Einstein product and its applications, refer to the works of \cite{behera2017further} and \cite{brazell2013solving}, which provide in-depth discussions on its mathematical formulation, properties, and practical implementations.
		\begin{definition}
            Let $\mathcal{A} \in \mathbb{R}^{I_1 \times \ldots \times I_L \times K_1 \times \ldots \times K_N}$ and $\mathcal{B} \in \mathbb{R}^{K_1 \times \ldots \times K_N \times J_1 \times \ldots \times J_M}$, The Einstein product of the tensors $\mathcal{A}$ and $\mathcal{B}$ yields a tensor of dimensions $\mathbb{R}^{I_1 \times \ldots \times I_L \times J_1 \times \ldots \times J_M}$. The elements of this tensor are given by the following expression:

$$
\left(\mathcal{A} *_N \mathcal{B}\right)_{i_1 \ldots i_L j_1 \ldots j_M} = \sum_{k_1, \ldots, k_N} a_{i_1 \ldots i_L k_1 \ldots k_N} b_{k_1 \ldots k_N j_1 \ldots j_M}
$$
		\end{definition}
\noindent			Here we comment that for a given tensor $\mathcal{A} \in \mathbb{R}^{I_1 \times \ldots \times I_N \times J_1 \times \ldots \times J_M}$, the transpose of $\mathcal{A}$, denoted by $\mathcal{A}^T$, is the tensor whose elements are $b_{i_1 \ldots i_M j_1 \ldots j_M}=a_{j_1 \ldots j_N i_1 \ldots i_M}$.\\
			The identity tensor $\mathcal{I}=\left[\mathcal{I}_{i_1, \ldots, i_N, j_1, \ldots, j_N}\right] \in \mathbb{R}^{I_1 \times \ldots \times I_N \times I_1 \times \ldots \times I_N}$ under the Einstein product is defined by the following entries
			$$
			\mathcal{I}_{i_1, \ldots, i_N, j_1, \ldots, j_N}= \begin{cases}1, & \text { if } i_k=j_k, \text { for } k\in\{1,2, \ldots, N\} \\ 0, & \text { otherwise. }\end{cases}$$
 The tensor $\mathcal{A} \in \mathbb{R}^{I_1  \times \cdots \times I_N\times I_1  \times \cdots \times I_N}$ is called the inverse of 
$\mathcal{X} \in \mathbb{R}^{I_1  \times \cdots \times I_N\times I_1 \times \cdots \times I_N}$
if it satisfies 
$\mathcal{A} *_{N} \mathcal{X}  = \mathcal{X}  *_N \mathcal{A} = \mathcal{I}.$ It is denoted as $\mathcal{A}^{-1}.$

			\begin{definition}
		(Tensor flattening). Let $\mathcal{A} \in \mathbb{R}^{I_1 \times \ldots \times I_N \times J_1 \times \ldots \times J_M}$ be a tensor of order $N+M$. The elements of the matrix $A \in \mathbb{R}^{I_1 \ldots I_N \times J_1 \ldots J_M}$ obtained by flattening the tensor $\mathcal{A}$ are given by
			$$
			A_{i, j}=\mathcal{A}_{i_1, \ldots, i_N, j_1, \ldots, j_M},
			$$
			where
			$$
			\begin{aligned}
				& i=i_1+\sum_{p=2}^N\left(i_p-1\right) \prod_{q=1}^{p-1} I_q, \\
				& j=j_1+\sum_{p=2}^M\left(j_p-1\right) \prod_{q=1}^{p-1} J_q .
			\end{aligned}
			$$
			\end{definition}	
            \noindent			Here the indices $i_p$ and $j_p$ belong to their domains, i.e., $1 \leq i_p \leq I_p$ for $1 \leq p \leq N$, and $1 \leq j_q \leq J_p$ for $1 \leq q \leq M$. We define $A=\operatorname{mat}(\mathcal{A})$ and $\mathcal{A}=\operatorname{mat}^{-1}(A)$.
\begin{proposition}

   Let $\mathcal{A} \in \mathbb{R}^{I_1 \times \ldots \times I_N \times J_1 \times \ldots \times J_M}$ and $\mathcal{B} \in \mathbb{R}^{J_1 \times \ldots \times J_M \times K_1 \times \ldots \times K_L}$ be tensors of orders $N+M$ and $M+L$, respectively. Then
			$$
			\operatorname{mat}\left(\mathcal{A} *_M \mathcal{B}\right)=\operatorname{mat}(\mathcal{A}) \cdot \operatorname{mat}(\mathcal{B}).
			$$

\end{proposition}

			\begin{definition} The scalar product between two tensors and its corresponding norm are defined as follows. 
   
\noindent
	The inner product of two same size tensors $\mathcal{X}, \mathcal{Y} \in \mathbb{R}^{I_1 \times \cdots \times I_N}$ is defined by
			$$
			\langle\mathcal{X}, \mathcal{Y}\rangle=\sum_{i_1=1}^{I_1} \sum_{i_2=1}^{I_2} \ldots \sum_{i_N=1}^{I_N} x_{i_1 i_2 \cdots i_N} y_{i_1 i_2 \cdots i_N} .
			$$
	
   \noindent
   The corresponding Frobenius norm $\mathcal{X}$ is given by
			$$
\|\mathcal{X}\|_F=\sqrt{\sum_{i_1 \ldots i_N} |x|^{2}_{i_1 \ldots i_N i_1 \ldots i_N}}.
			$$
		
		\end{definition}
        \begin{definition}(Kronecker product of tensors).
The left Kronecker product of $N^{th}$-order tensors, $\mathcal{A}\in\mathbb{R}^{I_1\times I_2\times...\times I_N}$\;and $\mathcal{B}\in\mathbb{R}^{J_1\times J_2\times...\times J_N}$ is a tensor $$\mathcal{C}=\mathcal{A}\otimes_{L} \mathcal{B}\in\mathbb{R}^{I_1J_1\times I_2J_2\times...\times I_NJ_N},$$ 
elementwise,
$$ c_{\overline{i_1j_1},...,\overline{i_Nj_N}}=a_{i_1,...,i_N}b_{j_1,...,j_N}.$$
\end{definition}
 \section{\textbf{Multiplex network }}
        Multiplex networks represent a powerful way to model systems that involve multiple types of interactions between entities.  such as social networks, communication systems, and biological networks. Unlike traditional networks that represent only a single relationship between nodes, multiplex networks capture multiple layers, each representing a distinct type of connection. However, this multilayered structure introduces complexity in both representation and analysis. A highly effective approach to representing multiplex networks is to use fourth-order tensors \cite{bergermann2022fast,bentbib2025tensor}. In this framework, each layer \(k \in \{1, 2, \dots, L\}\) of the multiplex network is represented by a non-negative adjacency matrix \(A^{(k)} \in \mathbb{R}^{N \times N}\), where \(N\) is the number of nodes, and \(L\) is the number of layers. Each matrix captures the strength or presence of connections between nodes for a particular type of interaction, such as friendship, professional relationships, or shared interests.\\
\noindent
The tensor representation allows for the simultaneous analysis of these interdependencies, making it a powerful tool for studying how interactions in one layer may influence others. This method provides a unified framework that simplifies the modeling and analysis of complex, multilayered systems, opening up new possibilities for exploring the dynamics of multiplex networks.
 \begin{figure}[h!]
        \centering
        \includegraphics[width=0.65\linewidth]{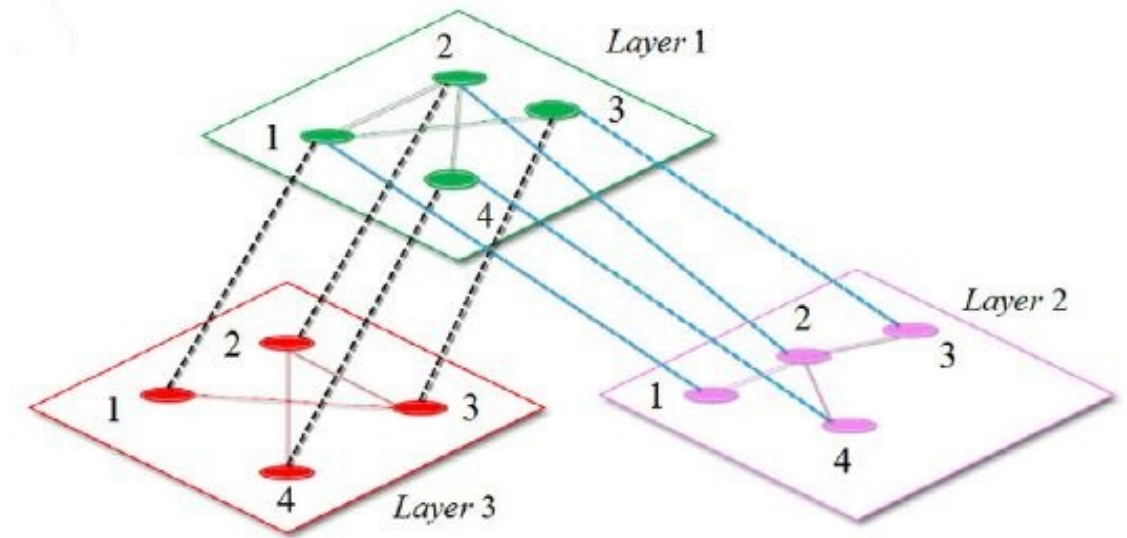}
        \caption{An example of a Multiplex network}
        \label{fig:enter-label}
    \end{figure}
\section{The Generalized analysis of random graph model approach and meta-algorithm for multiplex network structure}\label{sec3}
 For a multiplex network represented by \( G(V, E) \), where \( V \) is the set of vertices (nodes) and \( E \) is the set of edges, and with \( N \) representing the total number of nodes, we define the adjacency tensor \( \mathcal{A} \in \mathbb{R}^{N \times L \times N \times L} \). This adjacency tensor encapsulates the structure of the entire multiplex network, where \( N \) is the number of vertices (nodes) and \( L \) is the number of layers in the network. Each layer \( k \) in the multiplex network represents a distinct relationship or interaction between the nodes, and the tensor \( \mathcal{A} \) is constructed to store information about these interactions across all layers.

In this framework, the elements of the adjacency tensor, \( a_{i\alpha j\beta} \), represent the presence of an edge between node \( v_i^\alpha \) in layer \( \alpha \) and node \( v_j^\beta \) in layer \( \beta \). The value of \( a_{i\alpha j\beta} \) is greater than zero if and only if there exists an edge connecting node \( v_i^\alpha \) to node \( v_j^\beta \) in layer \( \alpha \), and is zero otherwise.

In this paper, our primary focus is on the random Kronecker graph model, which is a well established framework for generating synthetic network structures. A random Kronecker graph is defined as follows:

  \begin{definition}
   We define a graph \( G \) with \( N \) vertices as following a random Kronecker graph with initiator tensor \( \mathcal{P}_1 \) with probability:

\[
\mathcal{P}_1 = \{ P_{i\alpha j \beta} \}_{i, j \in \{ 1, \dots, N \} \atop \alpha, \beta \in \{ 1, \dots, L \}} \in \mathbb{R}^{N \times L \times N \times L},
\]
if the entries of its adjacency tensor \( \mathcal{A} \in \{0, 1\}^{N \times L \times N \times L} \) are (up to vertex correspondence via a permutation tensor \( \mathcal{\pi} \) of size \( N \)) independently drawn from a Bernoulli distribution with parameter \( \mathcal{P}_K \in \mathbb{R}^{N \times L \times N\times L} \). The probability tensor \( \mathcal{P}_K \) is the Kronecker power \( K \)-th of \( \mathcal{P}_1 \):
\[ \mathcal{P}_K = \mathcal{P}_{K-1} \otimes \mathcal{P}_1 = \mathcal{P}_1 \otimes \cdots \otimes \mathcal{P}_1 = \mathcal{P}_1^{\otimes K}, \] where \( K \in \mathbb{N} \). Thus, for \( i, j \in \{ 1, \dots, N \} \), \(\alpha, \beta \in \{ 1, \dots, L \}\).

\[
[\mathcal{A}]_{i\alpha j\beta} \sim \text{Bern}([\mathcal{P}_K]_{i\alpha j\beta}), \quad [\mathcal{A}]_{i\alpha j\beta} = [\mathcal{P}_K]_{i\alpha j\beta} + [\mathcal{Z}]_{i\alpha j\beta},
\]
where \( \mathcal{Z} \in \mathbb{R}^{N \times L \times N\times L} \) is a random tensor and it has independent entries of zero mean and variance. \\
\noindent \( \sum_{i,j,\alpha,\beta} [\mathcal{P}_K]_{ij\alpha\beta} \left( 1 - [\mathcal{P}_K]_{ij\alpha\beta} \right) \).
\end{definition}

  \begin{rmk}(Vertex matching) Consider a Kronecker graph \( G \) with \( N \) vertices as defined in Definition 7, where each vertex is given a unique label. The following relation holds:
\begin{equation}
    \mathcal{A} = \mathcal{\pi} *_{2}(\mathcal{P}_K + \mathcal{Z}) *_{2}\mathcal{\pi}^{-1},
\end{equation}
\noindent
for some permutation tensor \( \mathcal{\pi} \in \mathcal{P}_N \), where \( \mathcal{P}_N \) denotes the set of all permutation tensors.
  \end{rmk} 
\textbf{Assumption 1}: (Re-parameterization of the graph initiator) For fixed  \( L \), as \( N \to \infty \), the entry \( \mathcal{P}_{i\alpha j \beta} \) of the initiator \( \mathcal{P}_1 \) can be re-parameterized as:

\begin{equation}
\mathcal{P}_{i\alpha j \beta} = p + \frac{\mathcal{X}_{i\alpha j \beta}}{\sqrt{N}},
\end{equation}
\noindent
where \( p \in (0,1) \) and \( \mathcal{X} = \{\mathcal{X}_{i\alpha j \beta}\}_{i, j \in \{ 1, \dots, N \}, \alpha, \beta \in \{ 1, \dots, L \}} \in \mathbb{R}^{N \times L \times N \times L} \) with \( \|\mathcal{X}\|_{\max} = O(1) \).
\noindent
Next, we define the sequence of tensors \( \mathcal{S}_1, \dots, \mathcal{S}_K \) for \( \mathcal{S}_1 = \mathcal{X}/N \) and \( k \in \{2, \dots, K\} \) as:

\begin{equation}
\mathcal{S}_k = \frac{p^{k-1}}{N} *_{2} ( \mathcal{Y} *_{2} \mathcal{Y}^{\top} ) \otimes \mathcal{X} + p *_{2} \mathcal{S}_{k-1} \otimes (\mathcal{Y} *_{2} \mathcal{Y}^{\top}),
\end{equation}
\noindent
where \( \mathcal{Y} = \text{ones}(N^{k-1}, L^{k-1}, N^{k-1}, L^{k-1}) \).

\noindent
Finally, we show that the \( K \)-th Kronecker power \( \mathcal{P}_K = \mathbb{E}[\mathcal{A}] \) of \( \mathcal{P}_1 \) is closely related to \( \mathcal{S}_K \), which is directly associated with the graph parameters \( \mathcal{X} \), and has a low rank relative to its dimension \( N \). 
\begin{proposition}
{(Approximate small-rankness of \( \mathcal{P}_K \))} 

\begin{itemize}
    \item[(i)] \( \|\mathcal{P}_K - \mathcal{P}_K^{\text{lin}}\|_{\max} = O(N^{-1}) \) and \( \|\mathcal{P}_K - \mathcal{P}_K^{\text{lin}}\|_2 = O(1) \) for a linearized \( \mathcal{P}_K^{\text{lin}} \) defined as:
    \begin{equation}
        \mathcal{P}_K^{\text{lin}} = p^K \mathcal{W}*_{2} \mathcal{W}^{\top}  + \sqrt{N} \mathcal{S}_K,
    \end{equation}
  with \( \mathcal{S}_K \) in (3) for \( k = K \) so that \( \|\mathcal{S}_K\|_2 = O(1) \); and $\mathcal{W}=ones (N,L,N,L)$.
    
    \item[(ii)] \( \mathcal{S}_K \) is linear in (the entries of) \( \mathcal{X} \), in the sense that:
    \begin{equation}
        \mathcal{S}_K = \mathcal{\theta}*_2 \mathcal{X} \in \mathbb{R}^{N \times L \times N \times L},
    \end{equation}
    for known coefficients \( \mathcal{\theta} \in \mathbb{R}^{N \times L \times N \times L} \) (from binomial expansion). 
    
    \item[(iii)] \( \max(\operatorname{rank}(\mathcal{S}_K),  \operatorname{rank}(\mathcal{P}_K^{\text{lin}})) \leq (N-1)K + 1 \).
\end{itemize}    
\end{proposition}
\noindent	
\begin{theorem}
   {(Signal-plus-noise decomposition for \( \mathcal{A} \))} Given that Assumption 1 holds and  \( p^K \equiv p \in (0,1) \), the adjacency \( \mathcal{A} \) of a Kronecker graph in Definition 7 satisfies, for \( M \) large, \( \|\mathcal{A}\|_2 = O(\sqrt{N}) \) and
\begin{equation}
    \|\mathcal{A} - (\mathcal{\pi} *_2\mathcal{P}_K^{\text{lin}}*_2 \mathcal{\pi}^{-1} + \mathcal{Z})\|_2 = O(1),
\end{equation}
with
\begin{equation}
    \mathcal{\pi} *_2\mathcal{P}_K^{\text{lin}} *_2\mathcal{\pi}^{-1} = p^K \mathcal{W}*_{2} \mathcal{W}^{\top}  + \sqrt{N} \underbrace{\mathcal{\pi}*_2 \mathcal{S}_K *_2 \mathcal{\pi}^{-1}}_{\equiv \mathcal{S}_K^{\mathcal{\pi}}}.
\end{equation}
\end{theorem}
\begin{pf}
    According to Definition 7, we have \( \mathcal{A} = \mathcal{\pi} *_2 (\mathcal{P}_K + \mathcal{Z}) *_2\mathcal{\pi}^{-1} = \mathcal{\pi} *_2 \mathcal{P}_K^{\text{lin}}*_2 \mathcal{\pi}^{-1} + \mathcal{Z} + O_{\|\cdot\|_2}(1) \), where we take advantage of the fact that \( \|\mathcal{P}_K - \mathcal{P}_K^{\text{lin}}\|_2 = O(1) \) as shown in Proposition 3, and that the distribution of \( \mathcal{Z} \) remains unchanged after permutation by \( \pi \). This completes the proof of the Theorem.
\end{pf}
\begin{lemma}
    {(Consistent estimation of \( p \))} \textit{Under the notations and settings of Theorem 1, we have}
\begin{equation}
    \mathcal{W}^{\top} *_2\mathcal{A} *_2\mathcal{W} / N^2 - p^K \to 0, \quad \text{almost surely as } N \to \infty.
\end{equation} 
\end{lemma}
\begin{pf}
By Theorem 1, we have $\mathcal{A}=\mathbf{\pi} *_2\mathcal{P}_K *_2\mathbf{\pi}^{-1}+\mathcal{Z}$. First, we can note that $\mathcal{W}_N^{\top} *_2\mathcal{\pi}*_2 \mathcal{P}_K *_2\mathcal{\pi}^{-1} *_2\mathcal{W}_N=\mathcal{W}_N^{\top}*_2 \mathcal{P}_K*_2 \mathcal{W}_N=$ $p^K N^2+\sqrt{N} \mathcal{W}_N^{\top}*_2 \mathcal{S}_K *_2\mathcal{W}_N+O(N)=p^K N^2+O\left(N^{3 / 2}\right)$, where we used $\left\|\mathcal{S}_K\right\|_{\max }=O\left(N^{-1}\right)$ so that $\sqrt{N} \mathcal{W}_N^{\top}*_2 \mathcal{S}_K *_2\mathcal{W}_N=O\left(N^{3 / 2}\right)$. It then follows from the strong law of large numbers that $\frac{1}{N^2} \mathcal{W}_N^{\top}*_2 \mathcal{Z} *_2\mathcal{W}_N \rightarrow 0$ almost surely is as $N \rightarrow \infty$, and thus the conclusion.
\end{pf}
	
In plain words, Proposition 3 and Theorem 1 tell us that the adjacency tensor \( \mathcal{A} \) of a large Kronecker random graph can be decomposed, in a spectral norm sense, as the sum of some random zero-mean tensor \( \mathcal{Z} \) and (up to permutation by \( \mathcal{\pi} \)) the constant tensor \( p^K \mathcal{W} *_2\mathcal{W}^{\top} \) that can consistently estimate, according to Lemma 2) some deterministic signal tensor \( \mathcal{S}_K \) defined in (3).
\subsection{Kronecker denoising with shrinkage estimator}
In the context of graph analysis and tensor recovery, the adjacency tensor \( \mathcal{A} \) often contains both informative and noisy components. The primary goal is to denoise this random adjacency tensor in order to extract a low-rank tensor \( \mathcal{S}_K \) that captures the essential structure of the graph. This low-rank tensor, in turn, will help recover the underlying graph parameters \( \mathcal{X} \). To achieve this denoising process, we define the centered adjacency tensor \( \bar{\mathcal{A}} \), which is formulated as:

\begin{equation}
\bar{\mathcal{A}} \equiv \frac{1}{\sqrt{N}} \left( \mathcal{A} - \frac{\mathcal{W}^{\top} *_2 \mathcal{A} *_2 \mathcal{W}}{N^2} *_2 \mathcal{W} *_2 \mathcal{W}^{\top} \right),
\end{equation}

where \( N \) represents the number of nodes in the graph. We will demonstrate in the following result that the centered adjacency tensor \( \bar{\mathcal{A}} \) also conforms to a generalized signal-plus-noise model by subtracting the unwanted and non-informative constant tensor \( p^K \mathcal{W} *_2 \mathcal{W}^{\top} \) from \( \mathcal{A} \).

\begin{proposition}{(Signal-plus-noise decomposition for $\bar{\mathcal{A}}$)}  
Assume that \( \mathcal{X} \equiv \{\mathcal{X}_{i\alpha j \beta} \}_{i, j \in \{ 1, \dots, N \} \atop \alpha, \beta \in \{ 1, \dots, L \}} \) is "centered," meaning that: 
\[
\sum_{i, j \in \{ 1, \dots, N \} \atop \alpha, \beta \in \{ 1, \dots, L \}} \mathcal{X}_{i\alpha j \beta} = O(N^{-1/2}).
\]
\noindent
Then, the centered adjacency tensor \( \bar{\mathcal{A}} \), as defined in (9), satisfies the following:

\begin{equation}
    \|\bar{\mathcal{A}} - (\mathcal{S}_K^{\mathcal{\pi}} + \mathcal{Z}/\sqrt{N})\|_2 = O(N^{-1/2}),
    \label{eq:12}
\end{equation}
\noindent
\textit{with small-rank $\mathcal{S}_K^{\pi}$ defined in (7) and random tensor $\mathcal{Z}$.}
\end{proposition}
\begin{pf}
    From lemma 1, we have: 
\[
\begin{aligned}
\mathcal{W}_N^{\top} *_2 \mathcal{A} *_2 \mathcal{W}_N
&= \mathcal{W}_N^{\top} *_2 \mathcal{P}_K *_2 \mathcal{W}_N + O(N) \\
&= \mathcal{W}_N^{\top} *_2 \mathcal{P}_K^{\operatorname{lin}} *_2 \mathcal{W}_N + O(N) \\
&= p^K N^2 + O(N).
\end{aligned}
\]
Moreover, by (ii) of Proposition 3 and the assumption that \( \mathcal{W}_N^{\top} *_2\mathcal{X}*_2 \mathcal{W}_N = O(N^{-1/2}) \), we have 
\[
\begin{aligned}
\mathcal{W}_N^{\top} *_2 \mathcal{S}_K *_2 \mathcal{W}_N
&= \mathcal{W}_{N}^{\top} *_2 \mathcal{\theta} *_2 \operatorname{vec}(\mathcal{X}) \\
&= \frac{p^{K-1} K}{N^2} *_2 \mathcal{W}_{N}^{\top} *_2 \operatorname{vec}(\mathcal{X}) \\
&= O(\sqrt{N}),
\end{aligned}
\]
so that 
\[
\frac{1}{N^2} \mathcal{W}_N^{\top}*_2 \mathcal{A}*_2 \mathcal{W}_N = p^K + O(N^{-1}),
\]
and 
\[
\frac{\mathcal{W}_N^{\top}*_2 \mathcal{A} *_2\mathcal{W}_N}{N^2}*_2 \mathcal{W}_N *_2\mathcal{W}_N^{\top} = p^K \mathcal{W}_N *_2\mathcal{W}_N^{\top} + O_{\|\cdot\|_2}(1).
\]
This concludes the proof of Proposition 3.
\end{pf}
As a result of Proposition 3, in order to recover the desired signal tensor \( \mathcal{S}_K^{\mathcal{\pi}} \) from the noisy \( \bar{\mathcal{A}} \), we formulate the following optimization problem.

\begin{equation}
    \min_{\mathcal{S}_K^{\mathcal{\pi}} \in \mathbb{R}^{N \times L\times N\times L}} \| \bar{\mathcal{A}} - \mathcal{S}_K^{\mathcal{\pi}} \|,
    \quad \text{s.t.} \quad \text{rank}(\mathcal{S}_K^{\mathcal{\pi}}) \leq (N - 1)K + 1,
    \label{eq:13}
\end{equation}
    for a certain tensor norm \( \| \cdot \| \), which could be the Frobenius norm \( \| \cdot \|_F \), the spectral norm \( \| \cdot \|_2 \), or the nuclear norm \( \| \cdot \|_* \), with the rank constraint arising from (iii) of Proposition 3.
    
   \noindent	 
When the rank of \( \mathcal{S}_K^{\mathcal{\pi}} \) is known, the standard approach to solving \eqref{eq:13} is the hard thresholding singular value decomposition (SVD) estimator, defined as:
\begin{equation}
    \text{mat}(\hat{\mathcal{S}}_K^{\mathcal{\pi}}) = \sum_{i=1}^{\text{rank}(\mathcal{S}_K^{\mathcal{\pi}})} \hat{\sigma}_i \hat{u}_i \hat{v}_i^T,
    \label{eq:14}
\end{equation}
\noindent	
where \( (\hat{\sigma}_i, \hat{u}_i, \hat{v}_i) \) denotes the triple of singular values (arranged in decreasing order) along with the corresponding singular vectors left and right of \( \mathcal{A} \).

\noindent	
More generally, we introduce the \textit{shrinkage estimator}, which generalizes the hard thresholding SVD in \eqref{eq:14}, as follows:

\begin{equation}
    \text{mat}(\hat{\mathcal{S}}_K^{\mathcal{\pi}}) = \sum_{i=1}^{N} f(\hat{\sigma}_i) \hat{u}_i \hat{v}_i^T,
    \quad \text{with } f: \mathbb{R}_{+} \to \mathbb{R}_{+},
    \label{eq:15}
\end{equation}
\begin{theorem}
{Asymptotic Behavior of the Adjacency Spectrum}
\noindent
Let the empirical singular value distribution \( \mu_{A} \), defined as the normalized counting measure over the singular values \( \hat{\sigma}_i \) (listed in descending order) of the centered adjacency tensor \( \mathcal{A} \) from equation (9), where \( p^K \to p \in (0,1) \), converge weakly as \( N \to \infty \) to the limiting distribution given by the rescaled quarter-circle law:

\[
\mu(dx) = \frac{\sqrt{4p(1-p) - x^2}}{p(1-p)\pi} \cdot \mathbf{1}_{[0, 2\sqrt{p(1-p)}]}(x) dx, \tag{14}
\]
\noindent
where this convergence holds with high probability as \( N \to \infty \), known as the rescaled quarter-circle law. Furthermore, let \( \ell_i = \lim_{N \to \infty} \frac{\sigma_i(\text{mat}(\mathcal{S}_K^{\mathcal{\pi}}))}{\sqrt{p(1-p)}} \), where \( \sigma_i(\text{mat}(\mathcal{S}_K^{\mathcal{\pi}})) \) denotes the \( i \)-th largest singular value of \( \text{mat}(\mathcal{S}_K^{\mathcal{\pi}}) \), with associated left and right singular vectors \( u_i \) and \( v_i \). Then, the top singular values and their corresponding left and right singular vector triples \( (\hat{\sigma}_i, \hat{u}_i, \hat{v}_i) \) of \( A \) follow the phase transition behavior described below:

\[
\hat{\rho}_i \to 
\begin{cases} 
\sqrt{p(1-p)} \left( 2 + \ell_i^2 + \ell_i^{-2} \right), & \text{if } \ell_i > 1 \\
2\sqrt{p(1-p)}, & \text{if } \ell_i \leq 1 
\end{cases} \tag{15}
\]
\noindent
For \( 1 \leq i \leq \text{rank}(\text{mat}(\mathcal{S}_K^{\mathcal{\pi}})) \) and \( 1 \leq j \leq N \), the following limits for the singular vectors hold:

\[
(u_i^\top \hat{u}_j)^2 \to (1 - \ell_i^{-2}) \cdot \mathbf{1}_{\ell_i \geq 1} \cdot \mathbf{1}_{i = j}, \tag{16}
\]
\[
(v_i^\top \hat{v}_j)^2 \to (1 - \ell_i^{-2}) \cdot \mathbf{1}_{\ell_i \geq 1} \cdot \mathbf{1}_{i = j}.
\]
\end{theorem}
\begin{pf}
    The singular values \( \sigma_i(\tilde{A}) \) of the matrix \( \tilde{A} \in \mathbb{R}^{N \times N} \) correspond to the square roots of the eigenvalues \( \lambda_i(\tilde{A}\tilde{A}^\top) \), i.e., \( \sigma_i(\tilde{A}) = \sqrt{\lambda_i(\tilde{A}\tilde{A}^\top)} \). Therefore, to analyze the singular values and the associated eigenvectors, it is sufficient to evaluate the eigenvalues and eigenvectors of the positive semi-definite matrices \( \tilde{A}\tilde{A}^\top \) and \( \tilde{A}^\top\tilde{A} \). From Proposition 3, we know that as \( N \) grows large, \( \tilde{A} \) can be decomposed as the sum of a zero-mean random matrix \( Z/\sqrt{N} \) and a small-rank deterministic signal matrix \( S^{\Pi}_K \), i.e.,

\[
\tilde{A} = Z/\sqrt{N} + S^{\Pi}_K + \mathcal{O}_{\ell_2}(N^{-1/2}).
\]
\noindent
The asymptotic behavior of the eigenvalues and eigenvectors of sample covariance matrices, such as \( \tilde{A}\tilde{A}^\top \) or \( \tilde{A}^\top\tilde{A} \), follows directly from this decomposition, or the matrices \( \tilde{A}^\top \tilde{A} \) are quite standard in random matrix theory, but this only holds when the rank of the signal matrix \( S^{\Pi}_K \) is fixed with respect to its dimension \( N \to \infty \), as demonstrated in various studies (see, for example, \cite{bai2012sample,gavish2017optimal} and the references therein).
\noindent
In contrast, we consider the case where the rank of \( S^{\Pi}_K \) increases with the dimension \( N \), though at a slow rate such that \( \text{rank}(S^{\Pi}_K) \leq (N-1)K + 1 = (N-1) \log(N) + 1 = O(N) \), as shown in (iii) of Proposition 3. In this setting, it is sufficient to apply the deterministic equivalence results, such as those found in \cite{couillet2022random}, and observe that the resulting approximation errors scale as \( O(\log(N) N^{-1/4}) \) when \( Z \) has bounded and sub-Gaussian entries. This completes the proof of Theorem 2. 
\end{pf}
\begin{corollary}
    {(Shrinkage estimation of small-rank \( \mathcal{S}_K \))} We define the following shrinkage estimator:
\[
\text{mat}(\hat{\mathcal{S}}_K) = \sum_{i=1}^{N} f(\hat{\sigma}_i) \hat{u}_i \hat{v}_i^T,
\]
\noindent	
for \( f(t) = \sqrt{t^2 - 4p(1-p)} \cdot 1_{t > 2\sqrt{p(1-p)}} \) and \( (\hat{\sigma}_i, \hat{u}_i, \hat{v}_i) \) the triple of singular values and singular vectors of \( \bar{\mathcal{A}} \).\\
Suppose all singular values of \( \text{mat}(\mathcal{S}^\pi_K )\) that are greater than \( \sqrt{p(1-p)} \) are all distinct, one has:
\small
\[
\| \text{mat}(\mathcal{S}^\pi_K) - \text{mat}(\hat{\mathcal{S}}_K) \|_F^2 \to \sum_{i=1}^{\text{rank}(\text{mat}(\hat{\mathcal{S}}_K))} g(\sigma_i(\text{mat}(\mathcal{S}^\pi_K))) \quad \text{almost surely as} \quad N \to \infty,
\]
\noindent
with

\[
g(t) =
\begin{cases}
p(1-p) \left( 2 - p(1-p)t^{-2} \right), & t > \sqrt{p(1-p)}, \\
t^2, & t \leq \sqrt{p(1-p)}.
\end{cases}
\]
\end{corollary}
\begin{pf}
 We prove Corollary 1 by following the approach outlined in the proof of the theorem in \cite{gavish2017optimal}. By expanding the Frobenius norm, we obtain the following expression for \( \sigma_i \equiv \sigma_i(\mathcal{S}_K^\pi) \), the ordered singular values of \( \mathcal{S}_K^\pi \), and \( f(t) \). 
\[
\| \text{mat}(\mathcal{S}_K^\pi) - \text{mat}(\hat{\mathcal{S}}_K) \|_F^2 = \sum_{i=1}^r \left[ (\sigma_i)^2 + (f(\hat{\sigma}_i))^2 \right] - 2 \sum_{i,j=1}^r \sigma_i f(\hat{\sigma}_i) (u_i^{\top}u_j)(v_i^{\top}v_j) + O(1)
\]
\[
= \sum_{i=1}^r \left[ (\sigma_i)^2 - 2 \sigma_i f(\hat{\sigma}_i)(u_i^{\top}u_i)(v_i^{\top}v_i) + (f(\hat{\sigma}_i))^2 \right] + O(1),
\]
In the first equality, we use the fact that, by Theorem 2, there are at most \( r = \text{rank}(\text{mat}(\mathcal{S}_K) )\) singular values \( \hat{\sigma}_i \) of \( A \) greater than \( 2 \sqrt{p(1-p)} \). In the second and third lines, we apply the asymptotic singular vector characterization from Theorem 2. It can be demonstrated that the nonlinear shrinkage estimator \( \text{mat}(\hat{\mathcal{S}}_K) \) introduced in Corollary 1 minimizes the (asymptotic) Frobenius norm error among all estimators of the form \( \text{mat}(\hat{\mathcal{S}}_K) = \sum_{i=1}^N f(\hat{\sigma}_i) \hat{u}_i \hat{v}_i^{\top} \), where \( f: \mathbb{R}_{+} \to \mathbb{R}_{+} \).
\end{pf}
\begin{algorithm}[H]
\caption{Shrinkage estimator of $\hat{\mathcal{S}}_K$ to denoise $\mathcal{A}$}\label{alg:A1}
\begin{algorithmic}[1]
\Require The adjacency tensor $\mathcal{A}$ of a random Kronecker graph having $N$ nodes as in Definition 7.
\Ensure Shrinkage estimator $\hat{\mathcal{S}}_K$ of the (permuted) signal tensor $\mathcal{S}_K^{\mathcal{\pi}}$ defined in (7).
\State Compute the ``centered'' adjacency $\bar{\mathcal{A}}$ as in (9).\\
\State Estimate $p$ with $p = \frac{\mathcal{W}^{\top} *_2\mathcal{A} *_2\mathcal{W}^{\top}}{N^2}$ as in Lemma 1.\\
\State \Return $\text{mat}(\hat{\mathcal{S}}_K) = \sum_{i=1}^{(NL-1) \log (N) + 1} f(\hat{\sigma}_i)\hat{u}_i \hat{v}_i^T$, with $(\hat{\sigma}_i, \hat{u}_i, \hat{v}_i)$ the triple of singular values (in decreasing order) and singular vectors of $\bar{\mathcal{A}}$, for $f(t) = \sqrt{t^2 - 4p(1 - p)} \cdot \mathbf{1}_{t > 2\sqrt{p(1 - p)}}$.\\
\State \text{Then} \;$ \hat{\mathcal{S}}_K=reshape(\text{mat}((\hat{\mathcal{S}}_K), N,L,N,L)$
\end{algorithmic}
\end{algorithm}

\subsection{Kronecker solving via permuted linear regression}
For large values of \( N \) (the number of nodes in the network), the estimation of the low-rank tensor \( \hat{\mathcal{S}}_K \) that represents the essential structure of a Kronecker graph can be efficiently achieved using permutation-based techniques. We have \( \hat{\mathcal{S}}_K \approx \mathcal{\pi} *_2 \mathcal{S}_K *_2 \mathcal{\pi}^{-1} \), with

\begin{equation}
    \text{vec}(\hat{\mathcal{S}}_K) \simeq \text{vec}(\mathcal{\pi} *_2\mathcal{S}_K*_2 \mathcal{\pi}^{-1}) = (\mathcal{\pi} \otimes \mathcal{\pi})*_2 \mathcal{\theta}*_2 \text{vec}(\mathcal{X}).
    \label{eq:16}
\end{equation}

\noindent	
In order to recover both \( \text{vec}(\mathcal{X}) \) and \( \mathcal{\pi} \), we turn to the following optimization problem:

\begin{equation}
    (\hat{\mathcal{\pi}}, \text{vec}(\hat{\mathcal{X}})) = \arg\min_{\mathcal{\pi}, \text{vec}(\mathcal{X} )} \| (\mathcal{\pi} \otimes \mathcal{\pi})*_2 \mathcal{\theta}*_2 \text{vec}(\mathcal{X}) - \text{vec}(\hat{\mathcal{S}}_K) \|_2^2.
    \label{eq:17}
\end{equation}
Let \( \mathcal{\pi} \in \mathcal{P}_N \) denote the permutation tensor corresponding to the \textit{true} matching of the \( N \) vertices, and let \( d_H(\mathcal{\pi}, \mathcal{I}_N) \equiv | \{ i,j,\alpha,\beta \in \mathbb{N} : [\mathcal{\pi}]_{ij\alpha\beta} = 0 \} | \) represent the Hamming distance between \( \mathcal{\pi} \) and the identity tensor, which quantifies the number of mismatched vertices. The permuted linear regression problem in \eqref{eq:17} can then be expressed as:

\begin{equation}
    \min_{\mathcal{\pi} \in \mathcal{P}_N, \text{vec}(\mathcal{X})} \| \text{vec}(\hat{\mathcal{S}}_K )- (\mathcal{\pi} \otimes \mathcal{\pi})*_2 \mathcal{\theta}*_2 \text{vec}(\mathcal{X}) \|_2^2,
    \label{eq:18}
\end{equation}
\begin{equation*}
    \text{s.t.} \quad d_H(\mathcal{\pi}, \mathcal{I}_N) \leq s.
\end{equation*}
\\
This optimization problem is known to be NP-hard when \( s = O(N) \), except in the trivial case where \( N = 1 \).
\\
To solve the Kronecker inference problem efficiently, we consider the scenario where the permutation is \textit{sparse} (i.e., \( s \ll N \) in \eqref{eq:18}), and first relax the constraint in \eqref{eq:18} to \( d_H(\mathcal{\pi} \otimes \mathcal{\pi}, \mathcal{I}_N) \leq 2sN \). By introducing \( \text{vec}(\mathcal{D}) = (\mathcal{\pi} \otimes \mathcal{\pi} - \mathcal{I}_N) *_2 \mathcal{\theta} *_2 \text{vec}(\mathcal{X}) \), the problem in \eqref{eq:18} can be relaxed as:

\begin{equation}
    \min_{\text{vec}(\mathcal{X}), \text{vec}(\mathcal{D})} \| \text{vec}(\hat{\mathcal{S}}_K) - \mathcal{\theta} *_2\text{vec}(\mathcal{X}) - \text{vec}(\mathcal{D}) \|_2^2,
    \quad \text{s.t.} \quad \| \text{vec}(\mathcal{D} )\|_0 \leq 2sN,
    \label{eq:19}
\end{equation}
\noindent	
This remains non-convex because of the \( \ell_0 \) norm constraint. To address this, we explore the following two approaches:

\begin{itemize}
    \item[(i)] The iterative hard thresholding (IHT) method, which directly tackles the non-convex \( \ell_0 \)-norm constraint, using the hard thresholding operator \( H_s(\cdot) \) to set the tensor entries with small magnitudes to zero while keeping the larger ones unchanged; or
    \item[(ii)] Relaxing the non-convex problem in \eqref{eq:19} by replacing the \( \ell_0 \)-norm with the \( \ell_1 \)-norm, leading to the following Lagrangian formulation,
\end{itemize}

\begin{equation}
    \min_{\text{vec}(\mathcal{X}), \text{vec}(\mathcal{D})} \| \text{vec}(\hat{\mathcal{S}}_K) - \mathcal{\theta} *_2\text{vec}(\mathcal{X}) - \text{vec}(\mathcal{D}) \|_2^2 + \gamma \|\text{vec}( \mathcal{D} )\|_1,
    \label{eq:20}
\end{equation} 
    for a hyperparameter \( \gamma > 0 \) that balances the mean squared loss with the sparsity level in \( \mathcal{D} \).
    \begin{algorithm}[H]
\caption{Permuted Linear Regression to Solve for $\mathcal{X}$}\label{alg:A2}
\begin{algorithmic}[1]
\Require Estimated $\hat{\mathcal{S}}_K$ (from Algorithm 1), coefficient $\mathcal{\theta}$ and hyperparameter $\gamma$ for convex relaxation or step length $\eta$ and sparsity level $s$ for IHT.
\Ensure Estimation of graph parameter $\text{vec}(\hat{\mathcal{X}})$ by solving the permuted linear regression in (17).
\State \textbf{Initialize} $(\text{vec}(\hat{\mathcal{X}}), \text{vec}(\hat{\mathcal{D}}))$.
\State Estimate $p$ with $p = \frac{\mathcal{W}^{\top} *_2\mathcal{A} *_2\mathcal{W}^{\top}}{N^2}$ as in Lemma 1.
\State \Return $\text{mat}(\hat{\mathcal{S}}_K) = \sum_{i=1}^{(N-1) \log(N) + 1} f(\hat{\sigma}_i)*_2 \hat{u}_i *_2\hat{v}_i^T$, with $(\hat{\sigma}_i, \hat{u}_i, \hat{v}_i)$ the triple of singular values (in decreasing order) and singular vectors of $\bar{\mathcal{A}}$, for $f(t) = \sqrt{t^2 - 4p(1 - p)} \cdot \mathbf{1}_{t > 2\sqrt{p(1 - p)}}$.
\State \text{Then} $ \hat{\mathcal{S}}_K=reshape(\text{mat}((\hat{\mathcal{S}}_K), N,L,N,L)$.
\State \textbf{while} not converged
\\
\textbf{Option (I) IHT:}
\State $\text{vec}(\hat{\mathcal{Q}}) \leftarrow (1 - \eta) \text{vec}(\hat{\mathcal{D}}) + \eta (\text{vec}(\hat{\mathcal{S}}_K) - \mathcal{\theta} *_2\text{vec}(\hat{\mathcal{X}}))$
\State $\text{vec}(\hat{\mathcal{D}}) \leftarrow$ project $\text{vec}(\hat{\mathcal{Q}})$ into the set of sparse vectors via hard thresholding as $\text{vec}(\hat{\mathcal{D}}) = H_s(\text{vec}(\hat{\mathcal{Q}}))$
\\
\textbf{Option (II) Convex relaxation:}
\State $\text{vec}(\hat{\mathcal{D}}) \leftarrow \arg \min_{\text{vec}(\mathcal{D})} \| \text{vec}(\hat{\mathcal{S}}_K) - \mathcal{\theta} *_2\text{vec}(\mathcal{X}) - \text{vec}(\mathcal{D}) \|_2^2 + \gamma \| \text{vec}(\mathcal{D}) \|_1$ via soft thresholding
\State $\text{vec}(\hat{\mathcal{X}}) \leftarrow (\mathcal{\theta}^{\top} *_2 \mathcal{\theta})^{\dagger} *_2 \mathcal{\theta}^{\top} *_2 (\text{vec}(\hat{\mathcal{S}}_K) - \text{vec}(\hat{\mathcal{D}}))$.
\State \textbf{end while.}
\State \textbf{return} $\hat{\mathcal{X}}$
\end{algorithmic}
\end{algorithm}
\noindent	
Specifically, the signal tensor \( \mathcal{S}_K \) exhibits two key properties that make it crucial in the process of network completion and recovery:

\begin{itemize}
    \item[(i)] Small Rank and Denoising Ability: The signal tensor \( \mathcal{S}_K \) is characterized by having a small rank, which distinguishes it from the random tensor \( \mathcal{Z} \) as discussed in Proposition 3. This small rank property allows \( \mathcal{S}_K \) to be effectively "extracted" from noisy observations, such as the tensor \( \mathcal{A} \), through a denoising procedure. The denoising process aims to reduce the noise in the observed data, isolating the true signal tensor \( \mathcal{S}_K \) and thereby enhancing the accuracy of subsequent computations. This makes the signal tensor particularly valuable in practical scenarios where data is often corrupted by noise.
    
    \item[(ii)] Recovery of Tensor Data Using Perturbed Linear Regression: Another important feature of \( \mathcal{S}_K \) is that it is linear in the entries of \( \mathcal{X} \), with known coefficients \( \mathcal{\theta} \). This linearity allows for the application of a perturbed linear regression model to recover the desired tensor \( \mathcal{X} \) from the signal tensor \( \mathcal{S}_K \) (or its estimate). The use of linear regression in this context is significant because it provides a straightforward, computationally efficient method for recovering the original data, even in the presence of perturbations or noise. This approach facilitates the accurate reconstruction of network information from incomplete or noisy datasets, which is crucial for applications such as link prediction, community detection, and structural analysis.
\end{itemize}
\noindent
This leads to the development of a two-step meta-algorithm, consisting of a generalized denoising procedure followed by a solving step, as outlined in Algorithm 3. This meta-algorithm is designed to approximate the inference process for the parameters of the generalized random Kronecker graph, effectively handling both the noise reduction and parameter estimation tasks in sequence. The first step, denoising, aims to clean the observed data, removing unwanted noise, while the second step involves solving for the graph parameters, ensuring that the final results are as accurate and reliable as possible despite the initial noise in the data.
\begin{algorithm}[H]
\caption{Meta-algorithm: approximate inference of random Kronecker graph parameters}\label{alg:A3}
\begin{algorithmic}[1]
\Require The adjacency tensor \( \mathcal{A} \) of a random Kronecker graph of size \( N \) as in Definition 7.
\Ensure Estimates \( p \) and \( \hat{\mathcal{X}} \) of the graph parameters \( p \in \mathbb{R} \) and \( \mathcal{X} \in \mathbb{R}^{N \times L\times N \times L} \) in Assumption 1.
\State  Estimate \( p \) as \( p = \sqrt[K]{\mathcal{W}^{\top} *_2\mathcal{A} *_2\mathcal{W} / N^2} \) from Lemma 1.\\
 \textbf{Denoise} the adjacency tensor \( \mathcal{A} \) to get an estimate \( \hat{\mathcal{S}}_K \) of \( \mathcal{S}_K^{\mathcal{\pi}} \) defined in (7) with, e.g., the shrinkage estimator in Algorithm 1.\\
   \textbf{Solve} a permuted linear regression problem (see (15) above for detailed expression) to obtain \( (\hat{\mathcal{\pi}}, \text{vec}(\hat{\mathcal{X}})) \) from \( \hat{\mathcal{S}}_K \) via, e.g., the convex relaxation or the iterative hard thresholding approach in Algorithm 2.\\
   \textbf{return} \( p \) and \( \hat{\mathcal{X}} \).
\end{algorithmic}
\end{algorithm}
\section{Experiment Result}\label{sec5}
    In this section, we present examples of synthetic datasets to showcase the effectiveness of the previously described method. The calculations were performed using Matlab R2018b Intel(R) Core(TM) i5-8250U CPU @ 1.60GHz 1.80 GHz computer with 32 GB of RAM.
\\
\noindent
We analyse synthetic multiplex networks and compare the theoretical and empirical shrinkage estimation errors for the small-rank tensor \( \mathcal{S}_K \) as a function of the number of nodes \( N \). With an increase in \( N \), both the theoretical (red curve) and empirical (blue curve) errors decrease, indicating improved accuracy in the shrinkage estimation. The curves converge with larger 
\( N \), demonstrating that the theoretical model closely matches the synthetic data.
\begin{figure}[h!]
        \centering
        \includegraphics[width=0.5\linewidth]{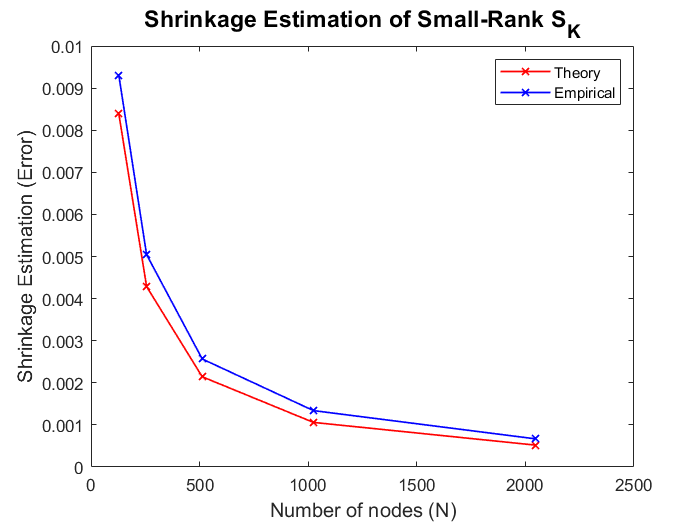}
        \caption{Shrinkage Estimation of Small-Rank tensor $\mathcal{S}_K$: Theory vs Empirical Results}
    \end{figure} 
    \newpage
    \noindent
    The figure below illustrates the operator norm error of the expression \( \mathcal{\bar{A}} - \left( \frac{\mathcal{Z}}{\sqrt{N}} + \mathcal{S}_K \right) \) as a function of \( N \), where \( N \) denotes the sample size. The plot displays the norm error against \( N \)  using a series of blue markers connected by a line. The norm error initially decreases rapidly as \( N \) increases, reflecting an improvement in the accuracy of the operator norm as the sample size grows. As \( N \) continues to increase, the error reduction becomes more gradual, leveling off at smaller values, indicating that the operator norm error approaches a stable minimum. The plot demonstrates the effectiveness of increasing \( N \) in reducing the error in the given model. 

    \begin{figure}[h!]
        \centering
        \includegraphics[width=0.5\linewidth]{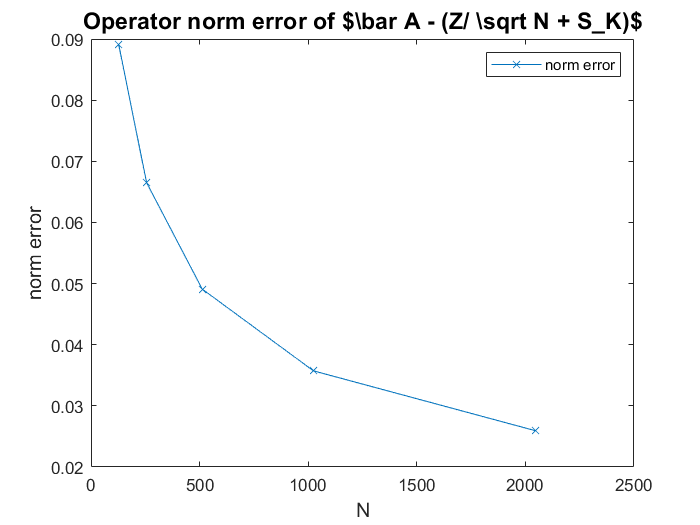}
        \caption{Operator norm error analysis of \( \bar{\mathcal{A}} - \left( \frac{\mathcal{Z}}{\sqrt{N}} + \mathcal{S}_K \right) \) with varying nodes sizes $N$}
    \end{figure} 
    \noindent
Furthermore, in our analysis, we applied the singular value decomposition (SVD) to the estimated matrix \( \text{mat}(\mathcal{S^{\pi}_K} )\). The first estimate of \( \text{mat}(\mathcal{S^{\pi}_K} )\)was obtained by applying hard thresholding (HS) on the SVD of the noisy centered adjacency matrix \( \tilde{A} \), as follows:

\[
\text{mat}(\mathcal{\hat{S}^{HS}_K}) = \sum_{i=1}^{rank(S^K)} \hat{\sigma}_i u_i \hat{v}_i^\top
\]
\noindent
where \( \hat{\sigma}_i \), \( u_i \), and \( \hat{v}_i \) represent the singular values (ordered in decreasing order) and the corresponding singular vectors of the matrix \( \tilde{A} \). Despite the simplicity of this estimator, it minimizes the difference in spectral norm \( \| \text{mat}(\mathcal{\hat{S}^{HS}_K}) - \text{mat}(\mathcal{\tilde{A})} \|_2 \) under the rank constraint \( rank(S_K) \), according to the Eckart-Young-Mirsky theorem. 
\begin{itemize}
    \item When the signal-to-noise ratio (SNR) \( \ell_i = \frac{\sigma_i(\text{mat}(\mathcal{S^{\pi}_K))}}{\sqrt{p(1-p)}} \) of \( S^{\Pi}_K \), as defined in Theorem 2, is below the phase transition threshold of 1, the corresponding estimated singular value \( \hat{\sigma}_i \) becomes independent of the true value \( \sigma_i(\text{mat}(\mathcal{S^{\pi}_K))} \). Furthermore, singular vectors asymptotically become orthogonal to the true vectors \( u_i \) and \( v_i \).

    \item Even when the SNR exceeds the threshold, the estimated singular values \( \hat{\sigma}_i \) remain different from the true values \( \sigma_i(\text{mat}(\mathcal{S^{\pi}_K))} \), and there exists a non-negligible angle between the estimated vectors \( \hat{u}_i \) and \( u_i \), as well as between \( \hat{v}_i \) and \( v_i \), unless the SNR \( \ell_i \) approaches infinity.
\end{itemize}
\noindent
In the figure below, we display the histogram of the singular values of \( \text{mat}(\mathcal{\tilde{A}}) / \sqrt{p(1 - p)} \) (in blue) compared to the limiting quarter-circle law spectrum and spikes (in red). The plot visually contrasts the singular value distribution from the matrix \( \tilde{A} \) with the theoretical limiting behavior, where the singular values deviate from the expected quarter-circle law, exhibiting spikes as expected from the noisy data. The analysis is performed with \( m = 2 \), \( K = 10 \), resulting in \( N = 4096 \), with \( p = 0.8 \) and \( \text{vec}(\text{mat}(\mathcal{X})) = [-5.5, 5.5, -1.5, 1.5]^\top \).
\begin{figure}[h!]
        \centering
        \includegraphics[width=0.5\linewidth]{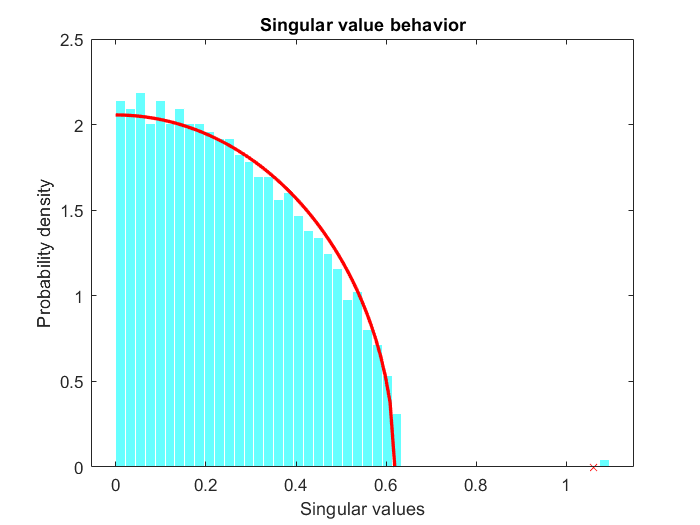}
         \caption{Histogram of singular values of ${A} \Big/ \sqrt{p(1 - p)}$}
    \end{figure} 

			\medskip
\section*{Conclusion}
In this paper, we introduce a highly efficient and scalable method for completing multidimensional networks using a generalized graph model. Our approach combines tensor decomposition techniques with advanced concepts from random tensor theory, providing a robust framework for addressing the challenges of multidimensional data. By utilizing tensor-based methods, we are able to decompose complex network structures into more manageable components, facilitating the accurate and efficient filling of missing or incomplete network data. This approach not only simplifies the network completion process but also significantly reduces computational costs, making it suitable for large-scale, real-world networks. The proposed framework ensures high accuracy in network analysis tasks.










\bibliographystyle{cas_model2_names}
\bibliography{ref}



\end{document}